\documentclass{article}

\usepackage{amsfonts}
\usepackage{latexsym}
\usepackage{amssymb}
\usepackage{amsmath}
\usepackage{amsthm}
\usepackage{verbatim}
\usepackage{array}
\usepackage{enumerate}
\usepackage{graphicx}           
\usepackage{graphics}           

\newtheorem{thm}{Theorem}[section]
\newtheorem{prop}[thm]{Proposition}
\newtheorem{cor}[thm]{Corollary}
\newtheorem{lem}[thm]{Lemma}

\newtheorem{defn}[thm]{Definition}
\newtheorem{remark}[thm]{Remark}
\newtheorem{example}{Example}
\newenvironment{rem}{\begin{remark}\rm}{\end{remark}}
\newenvironment{exa}{\begin{example}\rm}{\end{example}}

\newcommand{\QED}{\qed\medskip}

\newcommand{\nl}{\par\medskip\noindent}
\newcommand{\Proof}{{\nl\it Proof.\ }}
\newcommand{\Basis}{{\nl\textsc{Basis}.\ }}
\newcommand{\Hypothesis}{{\nl\textsc{Hypothesis}.\ }}
\newcommand{\Induction}{{\nl\textsc{Induction}.\ }}

\newcommand{\bt}{\mathbf{t}}
\newcommand{\bm}{\mathbf{m}}
\newcommand{\bx}{\mathbf{x}}

\newcommand{\initialSymbol}{\overline{\alpha}}
\newcommand{\familyT}{\mathcal{T}}
\newcommand\position[2]{\stackrel{\stackrel{#1}{\downarrow}}{#2}}

\newcommand{\N}{\mathbb N}
\newcommand{\Z}{\mathbb Z}
\newcommand{\Q}{\mathbb Q}

\newcommand{\Index}{\textsc{Index}}

\def\empt{\varepsilon}

\def\bbf{\mathbf{f}}

\begin{document}

\title{On the critical exponent of generalized Thue-Morse words\footnotemark[1]}

\author{Alexandre Blondin-Mass\'{e}\footnotemark[2],
    \thinspace  Sre\v{c}ko Brlek\footnotemark[3], \\
    \thinspace  Amy Glen\footnotemark[4],
    \thinspace S\'{e}bastien Labb\'{e}\footnotemark[5] \\[.3em]
{\small LaCIM, Universit\'e du Qu\'ebec \`a Montr\'eal,}\\
{\small C.P. 8888, succursale Centre-ville,
Montr\'eal, Qu\'ebec, CANADA, H3C 3P8}\\[.3em]
}

\date{Submitted: June 18, 2007; Accepted: October 15, 2007}

\maketitle

\footnotetext[2]{{\tt blondin$\_$masse.alexandre@courrier.uqam.ca}, NSERC (Canada) scholarship}
\footnotetext[3]{{\tt brlek@lacim.uqam.ca}, NSERC (Canada) grant}
\footnotetext[4]{{\tt amy.glen@gmail.com}, CRM-ISM-LaCIM postdoctoral fellowship (Montr\'eal, Canada)}
\footnotetext[5]{{\tt slabqc@gmail.com}, NSERC (Canada) scholarship}
\footnotetext[1]{Submitted to {\em Discrete Mathematics and Theoretical Computer Science}}

\begin{abstract}
    For certain generalized Thue-Morse words $\bt$, we compute the {\em critical exponent}, i.e., the supremum of the set of rational numbers that are exponents of powers in $\bt$, and determine exactly the occurrences of powers realizing it. \medskip
    
 \medskip

\end{abstract}

\section{Introduction}

    It is a well-known fact that the Norwegian mathematician Axel Thue (1863--1922) was the first to explicitly
    construct and study the combinatorial properties of an infinite {\em overlap-free word} over a 2-letter
    alphabet, obtained as the fixpoint of the morphism $\mu: \{a,b\}^*  
\rightarrow  \{a,b\}^* $
    defined by $\mu(a)= ab; \mu(b)=ba$:
    \[ \mu({\bm})={\bm} = abbabaabbaababba\cdots.\]
For a modern account of his papers, see Berstel~\cite{jB92axel}.
    Rediscovered by M.~Morse in 1921 in the study of symbolic dynamics, this
    overlap-free word is now called the \emph{Thue-Morse word}. 
This  ``ubiquitous'' sequence, already implicit in a memoir  by Prouhet~\cite{mP51memo} 
in 1851, appears in various fields, such as
    combinatorics on words, symbolic dynamics, differential geometry, number 
theory, and
    mathematical physics, as surveyed  by Allouche and Shallit~\cite{jAjS99theu}. That survey also mentions some generalizations of the 
Thue-Morse word, and recently other ones were considered in 
\cite{jAjS00sums, aF01over}.
  In particular,  the 
following result  was established in~\cite{jAjS00sums}.

    \begin{prop} \label{P:jAjS00sums} \emph{\cite{jAjS00sums}}
    Let $m \geq 1$, $b \geq 2$ be integers and let $s_b(n)$ denote the sum 
of the digits in the base $b$ representation of $n \in \N$. Then the 
infinite word
    $\bt_{b,m} := (s_b(n) \bmod{m})_{n\geq0}$ over the alphabet
    $\Sigma_m := \{0,1,\ldots, m-1\}$ is overlap-free if and only if $b
    \leq m$. \QED
    \end{prop}
It was also shown in \cite{jAjS00sums} that the word
    $\bt_{b,m}$ contains arbitrarily long squares, which extends a 
result previously established by Brlek~\cite{sB89enum} for $\bm$. Moreover, it was 
mentioned (\cite{jAjS00sums}, p.~8)  that,
    \begin{quote}
  {\em It would be interesting to determine the
    largest (fractional) power that occurs in the sequence $\bt_{b,m}$.
    For $b \leq m$, we already know that 2 is sharp.}
\end{quote}
We solve this problem here. Specifically, we study the family $\familyT$ of generalized 
Thue-Morse words consisting of the words $\bt_{b,m}$ as well as letter-renamings of them (see Section~\ref{S:definition}). For any $\bt \in \familyT$, we compute the {\em critical exponent} (i.e., the supremum of the set of rational numbers that are exponents of powers in $\bt$) and determine exactly the occurrences of powers realizing it, both in terms of $b$ and $m$. A noteworthy fact is that the critical exponents of generalized Thue-Morse words are always realized, which is not necessarily true in general (see Remarks~\ref{R:index} and \ref{R:krieger}). 

The next section contains all of the  basic terminology on words, borrowed 
mainly from Lothaire~\cite{mL83comb}, along with the generalized Thue-Morse 
words.  Section~\ref{S:preliminary}  contains the technical lemmas, 
establishing combinatorial properties used for proving our main results.   
The critical exponent (Theorem~\ref{T:critexp}) is then computed in Section 
\ref{S:critical}, and Section~\ref{S:occurrences} deals with the occurrences 
of factors realizing it (Theorem~\ref{T:genpos}).

Let us mention that the notions of `fractional  power' and `critical exponent' have received growing attention in recent times, especially in relation to {\em Sturmian} and {\em episturmian words}; see for instance \cite{jB99onth, vBcHlZ06init, dDdL02thei, aG05powe, jJgP01frac, jJgP02epis, dK06oncr, fMgP92repe, dV00stur}.

\section{Definitions and notation}

    Let $\Sigma$ denote a finite {\em alphabet}, i.e., a finite set of
    symbols called \emph{letters}. A {\em finite word} over $\Sigma$ is a
    finite sequence $w= w_{0}w_{1}\cdots w_{\ell -1}$, where each $w_i \in 
\Sigma$. We often write $w[i]$ for  $w_i$, when $i$ is a complicated formula. The
    \emph{length} of $w$ is $|w| = \ell$,
    and the {\em empty word} of length 0 is denoted by $\empt$. The set of all finite words over $\Sigma$ is denoted by $\Sigma^*$.

  A word $u$ is a \emph{factor} of $w$ if $w = pus$ for some words $p$, $s$. Moreover, such a word $u$ is said to be a {\em proper factor} of $w$ if $u\neq w$. We also say that $u$ is a {\em prefix} (resp.~{\em suffix}) of $w$ if $p = \empt$ (resp.~$s = \empt$). The number $|p|$ is called an {\em occurrence} of $u$ in $w$, i.e., $|p|$ is the beginning position of an appearance of $u$ in $w$.
    %
    The set of all factors of a word $w$ is denoted by $F(w)$.  A word $v$ is a {\em conjugate} of $w$ if there exists a word $u$ such that $uv = wu$.

An \emph{overlap} $w$ is a word of the form $w = auaua$ where $a$ is a letter and $u$ is a (possibly empty) word. For example, the English word {\em alfalfa} is an overlap
    and {\em banana} has the overlap {\em anana} as a suffix. A finite word is said 
to be \emph{overlap-free} if it does not contain an overlap as a factor.
Factors of the  form $uu$ are called {\em squares}.

The \emph{rational power} of a word $w$ is defined by
$w^r=w^{\lfloor r \rfloor}p$ where $r$ is a rational such that $r|w| \in \N$ and $p$ is the prefix of $w$ of length $(r-\lfloor r \rfloor)|w|$.
    For example, the word
    \[
    cab~\underline{abbacc}~\underline{abbacc}~\underline{abb}~bca
    \]
    contains a factor which is a $\frac{5}{2}$-power. 

      A (right) \emph{infinite word} (or simply a \emph{sequence}) $\bx$ over
    $\Sigma$ is a sequence indexed by $\N$ with values in $\Sigma$, i.e.,
    $\bx = x_0x_1x_2\cdots$ where each $x_i \in \Sigma$.  All of the 
terminology above naturally extends  to infinite words. 

An {\em ultimately periodic} infinite word can be written as $uv^\omega = uvvv\cdots$, for some $u$, $v \in \Sigma^*$, $v\ne \empt$. If $u = \empt$, then such a word is {\em periodic}. 
An infinite word that is not ultimately periodic is said to be {\em aperiodic}.

For any factor $w$ of an infinite word $\bx$, the
    \emph{index} of $w$ in $\bx$ is given by the number
    \[
    \Index(w) = \max\{r \in \Q  \mid w^r \in F(\bx)\}, 
    \]
    if such a number exists; otherwise, $w$ is said to have infinite
    index in $\bx$. The  {\it critical exponent} $E(\bx)$ of an infinite word is given by   \[
   {E(\bx) = \sup}\left\{\Index(w) \mid w \in F(\bx)\backslash\{\empt\}\right\}.
    \] 
It may be finite or infinite.   A  factor $w$ of $\bx$ is said to be a \emph{critical factor} if its index realizes the critical exponent of $\bx$, that is, when $\Index(w)=E(\bx)$.
    
\begin{rem} \label{R:index}
For any factor $w$ of an infinite word $\bx$, \Index$(w)$ is always realized when finite, whereas the critical exponent $E(\bx)$, even when finite, is not always realized; in particular, an infinite word may contain no critical factors. For example, the critical exponent of the well-known {\em Fibonacci word} $\bbf$ is  $E(\bbf) = 2 + \phi$, where $\phi = \frac{\sqrt{5}+1}{2}$ is the \emph{golden ratio} (see \cite{fMgP92repe}), but none of the factors of $\bbf$ realize $E(\bbf)$.
\end{rem}
    
    A \emph{morphism} is a function $\varphi : \Sigma^* \rightarrow \Sigma^*$ such 
that $\varphi(uv)=\varphi(u)\varphi(v)$ for all $u$, $v \in \Sigma^*$. For each letter $a \in \Sigma$, $\varphi (a)$ is called a \emph{block}.

\subsection*{Generalized Thue-Morse words} \label{S:definition}

    There exist many generalizations of the Thue-Morse word. Here, we 
introduce a morphism based
    formulation which is more convenient for our purposes. 

    \begin{defn} \label{D:gTMw}
        Let $b \geq 2$, $m \geq 1$ be integers,
        $\Sigma$ an alphabet of $m$ letters,
        $\sigma : \Sigma \rightarrow \Sigma$ a
        cyclic permutation and
        $\mu : \Sigma^* \rightarrow \Sigma^*$ the
        morphism given by
            $$\mu(\alpha) = \prod_{i=0}^{b-1} \sigma^i(\alpha)\
            = \sigma^0(\alpha)\sigma^1(\alpha)\sigma^2(\alpha)\cdots 
\sigma^{b-1}(\alpha).$$
        Then the generalized Thue-Morse word $\bt$,
        beginning with $\initialSymbol \in \Sigma$, is the infinite word 
given by
        $\bt = \mu^\omega(\initialSymbol)
               = {\lim}_{n\rightarrow\infty} \mu^n(\initialSymbol)$.
        \end{defn}

Hereafter, we study the family $\familyT$ of generalized Thue-Morse 
words $\bt$, for all $b  \geq 2$ , $m \geq 1$  and $\initialSymbol \in \Sigma$.

    \begin{rem} \label{R:nthletter}
        The $n$-th letter of $\bt$ is
        $\bt[n] :=\sigma^{s_b(n)}(\initialSymbol)$ where
        $s_b(n)$ denotes the sum of the digits in the base $b$
        representation of $n \in \N$.
    \end{rem}

  \begin{exa} \label{E:club}
  Let $b = 5$, $m = 3$,
    $\Sigma=\{\triangle,\diamondsuit,\heartsuit \}$ and
    $\sigma : \triangle \mapsto \diamondsuit \mapsto
             \heartsuit \mapsto \triangle $.
    This gives the following morphism :
    $\mu(\triangle)
    = \triangle\diamondsuit\heartsuit\triangle\diamondsuit$,
    $\mu(\diamondsuit)
    = \diamondsuit\heartsuit\triangle\diamondsuit\heartsuit$
    and
    $\mu(\heartsuit)
    = \heartsuit\triangle\diamondsuit\heartsuit\triangle$.
    By fixing $\initialSymbol=\diamondsuit$, we obtain the
    generalized Thue-Morse word
    $\bt = \mu^\omega(\diamondsuit) =
          \diamondsuit\heartsuit\triangle\diamondsuit\heartsuit
          \heartsuit\triangle\diamondsuit\heartsuit\triangle
          \triangle\diamondsuit\heartsuit\triangle\diamondsuit
          \diamondsuit\heartsuit\triangle\diamondsuit\heartsuit \cdots .$
      \end{exa}

  \begin{exa}
    With $\Sigma=\Z_m$, $\sigma(i) = (i + 1) \bmod m$ and
    $\initialSymbol=0$, Remark~\ref{R:nthletter} implies that
    \begin{equation} \label{E:AllSha}
        \bt[n]=\sigma^{s_b(n)}(0)=s_b(n) \bmod m,
    \end{equation}
    which is the case proposed by Allouche and
    Shallit~\cite{jAjS00sums}.
    \end{exa}

   \begin{exa}
     With $\Sigma=\Z_m$, $\sigma(i) = (i + c) \bmod m$  for some
    integer $c>0$ and $\initialSymbol=0$, we have 
    $\mu(0)=0\,\overline{c} \,\overline{2c} \, \cdots \, \overline{(b-1)c}$
    and
    \[\mu(i)=i \, \overline{(i+c)} \, \overline{(i+2c)} \,
    \cdots \, \overline{i+(b-1)c}=i^b+\mu(0)\]
    where $\overline{i} \equiv i \bmod m$.
    This is a subclass of the family of symmetric morphisms, defined by
    A. Frid~\cite{aF01over}. That paper also contains an extension of 
Proposition~\ref{P:jAjS00sums} to more general words.
    \end{exa}
    
    It is important to note that words in $\familyT$, up to letter renaming, are exactly those given in Equation (\ref{E:AllSha}). Our definition  avoids modular arithmetic on integers, and simplifies proofs by using the combinatorial properties  of $\familyT$ instead. For that purpose, we say
that a word $w = w_0 w_1 \cdots w_{\ell-1} \in F(\bt) $ is 
 \emph{$\sigma$-cyclic} if $w_{i}=\sigma(w_{i-1})$ for $1 \leq i \leq 
\ell-1$, or equivalently, if $w_i = \sigma^i(w_0)$
for $0 \leq i \leq \ell-1$.   As a consequence, blocks of $\bt$ are $\sigma$-cyclic.           

  \begin{exa}
    With $\Sigma=\{a,b\}$, $\sigma : a \mapsto b \mapsto a$  and
    $\initialSymbol=a$, we get $\mu(a) = ab$, $\mu(b)
    = ba$ and
    \[
    \bt = \mu^\omega(a) = abbabaabbaababba\cdots~
    \]
    which is the original Thue-Morse word $\bm$.
    \end{exa}

\begin{rem} \label{R:krieger}
Since $\bt \in \familyT$ is a fixpoint of a {\em uniform non-erasing} morphism, it follows immediately from Krieger's results in~\cite{dK06oncr} that the critical exponent of $\bt$ is either infinite (if $\bt$ is periodic) or rational. Moreover, since $\mu(\alpha)$ and $\mu(\beta)$ neither begin nor end with a common word when $\alpha \ne \beta$, it also follows from \cite{dK06oncr} that the critical exponent, when finite, is reached (i.e., $\bt$ contains critical factors).
\end{rem}
    
\section{Preliminary results}{\label{S:preliminary}}

    In this section, $\bt \in \familyT$ denotes an infinite
    generalized Thue-Morse word, as given in Definition~\ref{D:gTMw}.

    \begin{lem}{\label{L:period}} \emph{\cite{aF03arit, pMwM91digi}}
        The word $\bt$ is periodic if and only if $m \mid (b - 1)$.
        More precisely, if $\bt$ is periodic then
        $$\bt = \bigg[\prod_{i = 0}^{m - 1} \sigma^i(\initialSymbol)
    \bigg]^\omega. $$ \vspace{-1.1cm}
    \begin{flushright} \QED \end{flushright}
    \end{lem}


    \begin{lem}{\label{L:three}}
    Let $w$ be a $\sigma$-cyclic factor of $\bt$ of length $\ell$.
    If there exists an occurrence of $w$ overlapping three
    consecutive blocks, then $\bt$ is periodic.
    \end{lem}

    \Proof
        Suppose that there exists such an occurrence of $w$ in $\bt$. Since $w$ is 
        $\sigma$-cyclic and since every block is
        $\sigma$-cyclic, those three consecutive blocks,
        say $\beta_1$, $\beta_2$ and $\beta_3$, satisfy the fact
        that $\beta_1 \beta_2 \beta_3$ is $\sigma$-cyclic.

        Let $\alpha_1$, $\alpha_2$ and $\alpha_3$ be the first
        letters of $\beta_1$, $\beta_2$ and $\beta_3$. Then
        $\alpha_2 = \sigma^b(\alpha_1)$ and $\alpha_3 =
        \sigma^b(\alpha_2)$, and
        $\alpha_1 \alpha_2 \alpha_3 =
        \mu^{-1}(\beta_1 \beta_2 \beta_3)$
        are three consecutive letters occurring in $\bt$.
        This means, in particular, that either $\alpha_1
        \alpha_2$ or $\alpha_2 \alpha_3$ occur in the same block.
        Hence, either $\alpha_2 = \sigma(\alpha_1)$ or $\alpha_3 =
        \sigma(\alpha_2)$, which implies that 
        $\alpha_1 = \sigma^{b-1}(\alpha_1)$ or $\alpha_2 = \sigma^{b-1}
        (\alpha_2)$. Therefore, $m \mid (b - 1)$ and by Lemma
\ref{L:period},
        $\bt$ is periodic.
    \QED

    \begin{lem}{\label{L:cycl}}
        Let $w = w_0 w_1 \cdots w_{\ell-1}$ be a factor of $\bt$
        such that $b \nmid \ell$.
        \begin{enumerate}
            \item[{\em i)}]   \label{L:P:prefix}
                    If $wp$ occurs in $\bt$ where $p$ is any prefix
                    of $w$, then $p$ is $\sigma$-cyclic. 
            \item[{\em ii)}]    \label{L:P:power}
                    If $w^e$ occurs in $\bt$ for some rational $e>2$,
                    then $w^e$ is $\sigma$-cyclic.
        \end{enumerate}
    \end{lem} 
    \Proof
        Let $p$ be any prefix of length $\ell'$ of $w$ and assume $v = wp$
        occurs in $\bt$. Let $v_i$ be the $i$-th letter of $v$,
        with $0 \leq i \leq \ell + \ell' - 1$.
        For $0 \leq i \leq \ell' - 2$, either $v_i v_{i+1}$ or $v_{i+\ell} v_{i+\ell+1}$ is contained
        in the same block: if this was not the case, then $b$ would divide 
        $(i+\ell)-i=\ell$, contradicting our assumption. Since
        $v_i v_{i+1} = w_i w_{i+1} = v_{i+\ell} v_{i+\ell+1}$, we have 
        $w_{i+1} = \sigma(w_i)$ so that the first result follows.

     Since $w^e=www^{e-2}$, $w$ is $\sigma$-cyclic from i). Let $x$ be the conjugate of $w$ such that $w_0x=ww_0$. Then, $w^e=w_0xxx^{e-2-1/\ell}$ and from i) $x$ is $\sigma$-cyclic and so too is its factor $w_{\ell-1} w_{0}$. Hence, $w^e$ is $\sigma$-cyclic, which ends the second part.
    \QED

    \begin{lem}{\label{L:mult}}
        Suppose $\bt$ is aperiodic. If $w$ is a factor
        of $\bt$ of length $\ell \geq b$ with
        $b \nmid \ell$, then $\Index(w) \leq 2$.
    \end{lem}

    \Proof
        By contradiction, assume $v = www_0$ occurs in $\bt$, where
        $w_0$ is the first letter of $w$. Lemma~\ref{L:cycl} implies
        that $v$ is $\sigma$-cyclic. Moreover,
        $|v| = 2\ell + 1 > 2b$, which means that $v$ overlaps at least
        three consecutive blocks of $\bt$. Thus, by Lemma~\ref{L:three},
        $\bt$ is periodic which contradicts our assumption.
    \QED

\section{Critical exponent}{\label{S:critical}}

    In this section, we use exactly the same notation as previously.
    Before proving Theorem~\ref{T:critexp}, we need a few additional
    facts.

    The next three lemmas allow us to consider particular
    occurrences of factors of $\bt$.  We say that an occurrence $i$ of $w$ in $\bt$ is 
    \emph{synchronized} if $b \mid i$ and $b \mid (i + |w|)$.

    \begin{lem}{\label{L:leftright}}
        Suppose $w = w_0 w_1 \cdots w_{\ell-1}$ is a factor of $\bt$ such that $b \mid \ell$. 
        Let $p$ be a non-empty prefix of $w$ such that $bq+r$ is an occurrence of $w^np$ in $\bt$, where $q,r \in \N$, $0 \leq r < b$ and $n\geq 1$. Let $s$ be the possibly empty suffix of $w$ such that $w=ps$.
        Moreover, let $bq'+r'=bq+r+|w^np|$, where $q',r' \in \N$, $0 \leq r' < b$. Then,
        
        \begin{enumerate}
            \item[{\em i)}]
            $v=\bt[bq] \cdots  \bt[bq+r-1]$ is a suffix of $w$,
            \item[{\em ii)}]
            $u=\bt[bq'+r'] \cdots \bt[bq'+b-1]$ is a prefix of $sp$,
            \item[{\em iii)}]
            $bq$ is a synchronized occurrence of $vw^npu$ in $\bt$.
        \end{enumerate}
        
    \end{lem}

    \Proof
        Let $v = v_0 v_1 \cdots v_{r-1}=\bt[bq] \cdots \bt[bq+r-1]$ and let $s = s_0 s_1
        \cdots s_{r-1}$ be the suffix of $w$ of length $r$. We show that $v = s$.
        First note that $vw_0$ is contained in the block starting
        at position $bq$. Also, since $b \mid \ell$, 
        $sw_0$ is contained in one block. Both remarks imply that
        $vw_0$ and $sw_0$ are $\sigma$-cyclic. Hence,
        $$v_i = \sigma^{i-r}(w_0) = s_i \quad \textrm{for } 0 \leq i
        \leq r - 1,$$
        which gives the result. The proof of ii) is symmetric to the proof of i), and iii) follows from i) and ii).
    \QED

    \begin{lem}{\label{L:sync}}
        Let $w$ be a word of length $\ell$ such that $b \mid \ell$ and
        suppose $w^e$ occurs in $\bt$ for some rational $e > 1$.
        Then there exists a conjugate $x$ of $w$ and a rational
        $f \geq e$ such that both $x$ and $x^f$ have a synchronized
        occurrence in $\bt$.
    \end{lem}

    \Proof
        Let $p$ be the non-empty prefix of $w$ such that $w^e=w^n p$ where $n \in \N^{+}$ and let $s$ be the suffix of $w$ such that $w=ps$. Also, let $bq + r$ be an occurrence of $w^e$ in $\bt$, where $q,r \in \N$ and $0 \leq r < b$.
        By Lemma~\ref{L:leftright}, there exist a suffix $v$ of $w$ and a prefix $u$ of $sp$ such that $bq$ is a synchronized occurrence of $v w^e u$ in $\bt$.

        Now let $y$ be the prefix of $w$ such that $w=yv$, and define
        $x=vy$. We have $v w^e u = v (yv)^n p u = (vy)^n vpu = x^n
        vpu$. Moreover, $vpu$ is a prefix of $vpsps = vw^2 = v(yv)^2
        = x^2v$ so that $vw^eu = x^nvpu = x^f$ for some rational
        $f$. Finally,
        $$f = \frac{\ell n + |vpu|}{\ell} \geq
         \frac{\ell n+|p|}{\ell} = e,$$
        which ends the proof, since $x$ and $x^f$ both satisfy the required
        conditions.
       \QED

    The next lemma  deals with factors $w$ of $\bt$ of length $\ell$ not 
divisible by $b$.

    \begin{lem}{\label{L:basis}}
       Suppose $\bt$ is aperiodic and let $w$ be a factor
       of $\bt$ of length $\ell$ such that $b \nmid \ell$. Then
      \[
      \Index(w) \leq \left\{
        \begin{array}{ll}
            2b/m & \mbox{$~{\rm if~}\, b > m$}, \\
            2    & \mbox{$~{\rm if~}\,b \leq m$}.
        \end{array}
        \right .
      \]
    \end{lem}

    \Proof
        We distinguish three cases according to $\ell$, $b$ and $m$. \medskip
        
        \noindent \emph{Case $1:$ $\ell > b$}.
        By Lemma~\ref{L:mult}, we have $\Index(w) \leq 2$. \medskip

        \noindent \emph{Case $2:$ $\ell < b$ and $b > m$.}
        Let $w_0$ be the first letter of $w$ and $e = \Index(w)$.
        If $e > 2$, then by Lemma~\ref{L:cycl} 
        $w^e$ is $\sigma$-cyclic. In particular, $w_0 = \sigma^\ell(w_0)$,
        i.e., $m \mid \ell$ and $\ell \geq m$.
        However, by Lemma~\ref{L:three}, we know that $|w^e| \leq 2b$.
        Therefore,
        \begin{equation}\label{E:expon}
            e = \frac{|w^e|}{\ell} \leq \frac{2b}{\ell} \leq \frac{2b}{m}
        \end{equation}
        which ends this part. \medskip

       \noindent \emph{Case $3:$ $\ell<b$ and $b \leq m$.}
        As above, let $e = \Index(w)$ and suppose $e > 2$, which means $w^e$ 
is
        $\sigma$-cyclic and $m \mid \ell$. The latter statement is not
        possible since $1 \leq \ell < b \leq m$. We conclude in this
        case that $\Index(w) \leq 2$.
    \QED

    We are now ready to prove the main theorem of this section, which gives the critical exponent of $\bt$.     
        
      \begin{thm}{\label{T:critexp}}
        The critical exponent of $\bt$ is given by
\[ E(\bt) = \left \{
    \begin{array}{lll}
        \infty &\mbox{${\rm if}~\,m \mid (b-1)$,} \\
        2b/m &\mbox{${\rm if}~\,m \mid\hspace{-5pt}\not \hspace{5pt}(b-1)$ ~{\rm and}~ $b > m$,} \\
        2 &\mbox{${\rm if}~\,b \leq m$}.
        \end{array}
\right .
\]
    \end{thm}
    \Proof
        If $m \mid (b-1)$, then $E(\bt) = \infty$ since
        $\bt$ is periodic by Lemma~\ref{L:period}.
        Now suppose $m \nmid (b-1)$ and let
        \[
        E_{b,m} =  \left \{
    \begin{array}{lll}  
                           2b/m &\mbox{${\rm if~}\,b >m$,} \\
                           2 &\mbox{${\rm if~}\,b \leq m$}.
      \end{array}
      \right .
        \]
        First we show that $E(\bt) \geq E_{b,m}$.
        If $b \leq m$, it is easy to see that there is a square in the first 
two blocks
        of $\bt$, as noticed in~\cite{jAjS00sums}. On the other hand,
        if $b > m$, there exists a $2b/m$-power in $\bt$.
        Indeed, let $\beta_1$ and $\beta_2$ be the blocks starting
        at positions $b^m-b$ and $b^m$ respectively.
        From Remark~\ref{R:nthletter}, the last letter of $\beta_1$ is
        $\bt[b^m-1]=\sigma^{s_b(b^m-1)}(\initialSymbol)
                       =\sigma^{m(b-1)}(\initialSymbol)=\initialSymbol$
        and the first letter of $\beta_2$ is
        $\bt[b^m]=\sigma^{s_b(b^m)}(\initialSymbol)
                    =\sigma(\initialSymbol)$.
        Therefore, the whole factor $\beta_1\beta_2$ of length $2b$
        is $\sigma$-cyclic, and hence
        $\beta_1\beta_2=w^{2b/m}$ where $w$ is the prefix
        of length $m$ of $\beta_1$.
        Thus $E(\bt) \geq E_{b,m}$.

        We now prove that $E(\bt) = E_{b,m}$ by showing that $\Index(w) \leq
        E_{b,m}$ for any factor $w$ of $\bt$.
        Suppose $w$ is a factor of $\bt$ of length $\ell = b^iN$, where $b 
\nmid
        N$ for some $i, N \in \N$. The proof proceeds by induction on $i$.

        \Basis If $i = 0$ then $b \nmid \ell$. This case is proved in
        Lemma~\ref{L:basis}.

        \Hypothesis We assume $\Index(w) \leq E_{b,m}$ for all
        factors $w$ of $\bt$ of length $\ell = b^iN$.

        \Induction
        Let $w$ be a factor of $\bt$ of length $\ell = b^{i+1}N$. Assume $e = \Index(w)>1$.
        Since $b \mid \ell$, from Lemma~\ref{L:sync}
        we know that there exists a factor $x$ of length $\ell$ and a
        rational $f \geq e$
        such that both $x$ and $x^f$ are synchronized with the blocks.
        Let $p$ be the proper prefix of $x$ such that $x^f=x^{n}p$
        where $n \in \N^{+}$. Then,
        $\mu^{-1}(x^f)=\mu^{-1}(x^{n}p)=
        \left(\mu^{-1}(x)\right)^{n}\mu^{-1}(p)$
        is a factor of $\bt$. Since $\mu^{-1}(p)$ is
        a prefix of $\mu^{-1}(x)$, we obtain
        \begin{equation} \label{E:preimage}
            \Index\left(\mu^{-1}(x)\right)
            \geq n+\frac{|\mu^{-1}(p)|}{|\mu^{-1}(x)|}
            =n+\frac{|p|/b}{|x|/b}
            =n+\frac{|p|}{|x|}=\Index(x).
        \end{equation}
        By the induction hypothesis, we have 
            $$\Index(w) \leq f \leq \Index(x) \leq
\Index\left(\mu^{-1}(x)\right)
            \leq E_{b,m}.$$ \QED

\section{Occurrences of critical factors} {\label{S:occurrences}}

    In this section, we assume that $\bt$ is \emph{aperiodic},
    i.e., $m \nmid (b - 1)$. We say that $b>m$ is the \emph{overlap case}
    and that $b \leq m$ is the \emph{square case} and     
    denote by $e$ the critical exponent of $\bt$. 
    
    Here, we describe the \emph{occurrences} and the
    \emph{lengths} of the critical factors of $\bt$.

    \begin{lem}{\label{L:image}}
        Let $w$ be a critical factor of $\bt$ of length $\ell$.
        Then the following properties hold.
        \begin{enumerate}
        \item[{\em i)}]   $\mu (w)$ is a critical factor of $\bt$.
        \item[{\em ii)}]   If $b \mid \ell$, then $\mu^{-1} (w)$ is a critical
                factor of $\bt$.
        \end{enumerate}
    \end{lem}

    \Proof
        Property i) is trivial. For ii), suppose $b \mid \ell$. If
        $w^{e}$ is not synchronized with the blocks, then
        Lemma~\ref{L:leftright} contradicts
        the maximality of $e$. Therefore, both $w$ and $w^{e}$ are
        synchronized and their preimages under $\mu$ are well-defined.
        Finally, from inequality (\ref{E:preimage}),
        we obtain the index of the preimage.~\QED

    In view of Lemma~\ref{L:image}, it is enough to consider only
    the case $b \nmid \ell$ when describing the occurences of critical
    factors of length $\ell$ in $\bt$.

    \begin{lem}{\label{L:centered}}
        Let $w$ be a critical factor of $\bt$ of length $\ell$
        such that $b \nmid \ell$.
        \begin{enumerate}
            \item[{\em i)}]   In the overlap case, $\ell = m$ and $w^{e} = \beta_1
                    \beta_2$, where $\beta_1$, $\beta_2$ are two
                    consecutive blocks.
            \item[{\em ii)}]   In the square case, write $w^2 = w^{(1)}w^{(2)}$.
                    Then $w^{(2)}$ occurs at the beginning of a
                    block.
        \end{enumerate}
    \end{lem}

    \Proof
        \emph{Overlap case} ($b>m$). We consider separately the cases
        $\ell < b$ and $\ell > b$. If \emph{$\ell < b$}, then
        we have $\ell=m$ and $|w^{e}|=2b$ by inequality (\ref{E:expon}).
        It follows that $w^e$ is composed of two consecutive blocks.
        Otherwise, if $w^e$ overlaps three consecutive blocks,
        then from Lemma~\ref{L:three} $\bt$ is
        periodic, a contradiction.

        On the other hand, if \emph{$\ell>b$}, then from the last section 
        we know that $\Index(w)=2b/m>2$.
        This contradicts Lemma~\ref{L:mult}; hence there is
        no such critical factor. \medskip

        \noindent \emph{Square case} ($b \leq m$). Again, we distinguish two cases. First assume \emph{$\ell<b$}. Note that
        each block of $\bt$ contains distinct letters.
        Also, $w^2$ is a factor of $\bt$ and, from Lemma~\ref{L:cycl},
        $w$ is $\sigma$-cyclic.
        If $\alpha_{i}$ denotes the
        $i$-th letter of $w^2$ with  $0 \leq i < 2\ell$,
        then $\alpha_{i}=\alpha_{i+\ell}$ for all $0 \leq i < \ell$.
        This implies that
        $\alpha_{0}$ and $\alpha_{\ell}$ are not in the same block
        so that there exist two consecutive letters $\alpha_{i-1}$
        and $\alpha_{i}$, where $1 \leq i \leq \ell$, which are not
        contained in the same block. Suppose $1 \leq i < \ell$.
        Then, $\alpha_{i}$ and $\alpha_{i+\ell}$ belong
        to the same block since $\ell<b$. This is a contradiction
        because $\alpha_{i}=\alpha_{i+\ell}$.
        We conclude that the pair $\alpha_{\ell-1}$ and $\alpha_{\ell}$
        is not contained in the same block. In other words, $w^{(2)}$
        occurs at the beginning of a block.

        Lastly, consider the case \emph{$\ell>b$}.
        From Lemma~\ref{L:cycl}, $w$ is $\sigma$-cyclic.
        Suppose $w^{(2)}$ does not occur at the beginning of a
        block. Then the pair $w_{\ell-1}w_0$
        formed by the last and the first letter of $w$ is in the same
        block and is $\sigma$-cyclic. Hence, the whole factor $w^2$
        of length $2\ell>2b$ is $\sigma$-cyclic. This factor overlaps
        three consecutive blocks so $\bt$ is periodic from
        Lemma~\ref{L:three}, a contradiction.
    \QED

    It has already been noticed that squares of certain factors
    of length $3$ appear in the Thue-Morse word (see~\cite{sB89enum} for example). The following lemma proves the
    uniqueness of this fact.

    \begin{lem}{\label{L:only}}
        The Thue-Morse word ${\bm}$ is the unique word in $\familyT$ 
containing
        a critical factor $w$ of length $\ell > b$ such that $b \nmid \ell$.
        Moreover, $\ell=3$.
    \end{lem}

    \Proof
        Suppose that $w$ is critical factor of $\bt$ of length $\ell > b$
        such that $b \nmid \ell$. In the proof of Lemma
       ~\ref{L:centered}, we saw that in the overlap case
        there is no such critical factor.
        In the square case, we know also from Lemma~\ref{L:centered}
        that if $w^2 = w^{(1)} w^{(2)}$ then $w^{(2)}$ occurs at
        the beginning of a block. Moreover, Lemma~\ref{L:cycl} implies
        that $w$ is $\sigma$-cyclic. Since $\ell>2b$
        implies that $w$ overlaps three consecutive blocks,
        it follows from Lemma~\ref{L:three} that $b<\ell<2b$.

        Hence $w^2$ overlaps exactly four
        consecutive blocks, say $\beta_1$, $\beta_2$, $\beta_3$ and
        $\beta_4$. Let $\alpha_i$ be the first letter of the block
        $\beta_i$, $i \in \{1,2,3,4\}$, and let $r=2b-\ell$ be the
        distance between $\alpha_1$ and the first letter of $w^2$.
        Then $\alpha_3=w_0=\sigma^{r}(\alpha_1)$. Also, since $w$
        is $\sigma$-cyclic, we have $\alpha_2=\sigma^{b}(\alpha_1)$
        and  $\alpha_4=\sigma^{b}(\alpha_3)=\sigma^{b+r}(\alpha_1)$.
        Moreover, it follows that   
        $\alpha_1\alpha_2\alpha_3\alpha_4 =
        \mu^{-1}(\beta_1\beta_2\beta_3\beta_4)$
        is a factor of $\bt$. Assume $\alpha_1\alpha_2$
        is in the same block, then $\alpha_2=\sigma(\alpha_1)$; but then, since $\alpha_2=\sigma^{b}(\alpha_1)$, we get that $\alpha_1=\sigma^{b-1}\alpha_1$, and so $m \mid (b-1)$,
        which is in contradiction with $b \leq m$. Assuming $\alpha_3\alpha_4$ is in the same block gives the same contradiction. Hence, neither $\alpha_1\alpha_2$ nor
        $\alpha_3\alpha_4$ are in the same block. We conclude from the last
        observation that $b=2$ since
        $\alpha_2\alpha_3 = \sigma^{b}(\alpha_1)\sigma^{r}(\alpha_1)$
        must form a block. Therefore $\ell=3$ and $r=1$. Also, $m$ divides  
        $b+1-r=b$ and $m=b=2$ because the block $\alpha_2\alpha_3$
        is $\sigma$-cyclic, which ends the proof.

        Note that $\alpha_2=\alpha_1$,
        $\alpha_4=\alpha_3$, $\alpha_3=\sigma(\alpha_1)$ and $w =
        \sigma(\alpha_1) \alpha_1 \sigma(\alpha_1)$, which has already
        been noticed in~\cite{sB89enum}.
    \QED

    We now prove the following lemma which is a more general result
    than Lemma 5 in~\cite{jAjS00sums}.

    \begin{lem}{\label{L:sum}}
        Let $k, N \in \N^{+}$ be such that $b \nmid k$ and
        $1 \leq N < b$. Then,
        $$s_b(kb^q-N) - s_b(kb^q) = q(b-1) - N \quad \mbox{\rm for any $q \in
        \N^{+}$}.$$
    \end{lem}
    \Proof By direct computation,
    \begin{eqnarray*}
        s_b(kb^q-N)
              &=& s_b\left((k-1)b^q+(b^q-b)+(b-N)\right)\\
              &=& s_b\left[(k-1)b^q+(b-1)\sum_{i=1}^{q-1}b^{i}+\
                  (b-N)\right]\\
              &=& s_b(k-1) + (q-1)(b-1) + (b-N)\\
              &=& s_b(k) - 1 + q(b-1) - N + 1\\
              &=& s_b(kb^q) + q(b-1) - N, 
                  \end{eqnarray*}
                 as desired. \QED

    From {\em Bezout's Identity}, we know that for any $x, y \in \Z$ there exist
    $s, t \in \Z$ such that $\gcd(x,y) = sx+ty$ where $s$ can be chosen
    positive. Let $g\in \Z$ and 
    $$S^{g}_{x,y}=\{s \in \N^{+} \mid g=sx+ty, t \in \Z \}.$$
    The set $S^{g}_{x,y}$ is non-empty when $\gcd(x,y)\mid g$.
    Moreover, let us define the following three sets:
    \begin{eqnarray*}
        A     &=& \left\{ \left. kb^q-b \enskip \right| \
                k \in \N^{+}, b \nmid k\
                \mbox{and } q \in S^{m}_{b-1,m} \right\},\\
        B_{N} &=& \left\{ \left. kb^q-N \enskip \right| \
                k \in \N^{+}, b \nmid k\
                \mbox{and } q \in S^{N}_{b-1,m} \right\},\\
        C &=& ( 8\cdot B_1 +3) \cup ( 8\cdot B_1+7). \\
    \end{eqnarray*}
    
    We are now ready to state the main theorem of this section, which gives the set of occurrences where critical factors in an aperiodic generalized Thue-Morse word $\bt$ realize the critical exponent $e = E(\bt)$.
    
    
        \begin{thm}{\label{T:genpos}}
        If $w$ is a critical factor
    of $\bt$ of length $\ell=Nb^i$ such that $b \nmid N$,
    then the set of occurrences of $w^e$
    in $\bt$ is
    \[
 \left \{
    \begin{array}{lll}
        b^iA               & \mbox{${\rm if}~\,b > m$,} \\
        b^iB_{N}           & \mbox{${\rm if}~\,b \leq m$ ~{\rm and}~ $(b,m) \neq (2,2)$,}  \\
        b^i(B_{1} \cup C)  & \mbox{${\rm if}~\,b=m=2$.}\\
    \end{array}
    \right .
    \]
    \end{thm}

The proof of  Theorem~\ref{T:genpos} follows easily from Lemma~\ref{L:image} and the next  lemma.

    \begin{lem}{\label{L:pos}}
        Let $w$ be a critical factor of $\bt$ of length $\ell$
        such that $b \nmid \ell$.
        Then, the set of occurrences of $w^e$ in $\bt$ is
        \[
       \left\{ \begin{array}{ll} 
            A            &\mbox{${\rm if}~\,b > m$,} \\
            B_{\ell}     &\mbox{${\rm if}~\,b \leq m$ ~{\rm and}~ $(b,m) \neq (2,2)$,}  \\
            B_{1} \cup C &\mbox{${\rm if}~\,b=m=2$.}\\
        \end{array} \right.
        \]
    \end{lem}
    \Proof
       \emph{Overlap case} ($b > m$). From Lemma
       ~\ref{L:centered}, we know that $w^{e}=\beta_1\beta_2$, where 
        $\beta_1$ and $\beta_2$ are blocks. Suppose $p$ is an occurence of $w^e$
        and let $kb^q$ be the starting position of
        $\beta_2$ where $b \nmid k$ and $k,q \in \N^{+}$.
        Then $p=kb^q-b$. Now, let $\gamma_1=\bt[kb^q-1]$ be the last letter of $\beta_1$
        and $\alpha_2=\bt[kb^q]$ be the first letter of $\beta_2$.
        Since $w^{e}=\beta_1\beta_2$ is $\sigma$-cyclic from Lemma~\ref{L:cycl},
        we have
        $$\sigma(\sigma^{s_b(kb^q-1)}(\initialSymbol))
           = \sigma(\gamma_1)=\alpha_2 = \sigma^{s_b(kb^q)}(\initialSymbol).$$
        That is,
        \begin{eqnarray*}
            \initialSymbol
                &=& \sigma^{s_b(kb^q-1)-s_b(kb^q)+1}(\initialSymbol) \\
                &=& \sigma^{q(b-1) - 1 +1}(\initialSymbol)\\
                &=& \sigma^{q(b-1)}(\initialSymbol),
        \end{eqnarray*}
        and hence $m \mid q(b-1)$. Therefore, 
        there exists $t \in \Z$ such that $mt=q(b-1)$. In particular, we have   
        $m=q(b-1)-(t-1)m$ where $q \in S^{m}_{b-1,m}$
        which ends the first part. \medskip

        \noindent\emph{Square case} ($b \leq m$). If we write $w^2 = w^{(1)}w^{(2)}$, then from Lemma~\ref{L:centered} we know that $w^{(2)}$ occurs at the beginning of a block. We distinguish  two cases. 
        
        First consider the case $\ell<b$.   
        Suppose $p$ is an occurrence of $w^2$ and
        let $kb^q$ be the starting position of $w^{(2)}$ (i.e., the second
        block) where $b \nmid k$ and $k,q \in \N^{+}$. Then
        $p=kb^q-\ell$, and hence 
        $$\sigma^{s_b(kb^q-\ell)}(\initialSymbol)
           = \sigma^{s_b(kb^q)}(\initialSymbol),$$
        that is 
        \begin{eqnarray*}
            \initialSymbol
                &=& \sigma^{s_b(kb^q-\ell)-s_b(kb^q)}(\initialSymbol) \\
                &=& \sigma^{q(b-1) - \ell}(\initialSymbol).
        \end{eqnarray*}
        The last equation holds if and only if $m$ divides $q(b-1) - \ell$, 
        in which case there exists $t \in \Z$
        such that $\ell=q(b-1)-tm$, that is $q \in S^{\ell}_{b-1,m}$.

        Now consider $\ell>b$.
        From Lemma~\ref{L:only}, we know that such critical factors
        occur only in the Thue-Morse word $\bm$ and that
        $w^2 \in \{ abaaba, babbab \}$. By the fixpoint property,
        \[\bm = \mu^3(\bm) = \mu^3(a) \mu^3(b) \mu^3(b) \mu^3(a) \mu^3(b) 
\cdots,\] so that
       it  can be factorized into blocks of length $8$. Since
        $|w^2| = 6$,  
        $w^2$ is a factor of two consecutive blocks of length $8$.
        Those two blocks are either $\mu^3(\alpha \beta)$ or
        $\mu^3(\alpha \alpha)$, where $\alpha,
        \beta \in \{a,b\}$, $\alpha \neq \beta$. We observe that
        $$\mu^3(\alpha \beta) = \alpha \beta \beta \alpha \beta \alpha
            \alpha \beta \beta \alpha \alpha \beta \alpha \beta \beta 
\alpha$$
        contains no square of length $6$. On the other hand, $w^2$ occurs in
            $$\mu^3(\alpha \alpha) = \alpha \beta \beta \underline{\alpha} 
\beta \alpha
            \alpha \underline{\beta} \alpha \beta \beta \alpha \beta \alpha
            \alpha \beta$$
        at position $3$ or $7$.
        Thus, the set of occurrences of squares of length $6$ is 
exactly
        $(8\cdot B_1+3)\cup(8\cdot B_1+7)$, since $B_1$ enumerates
        the occurrences of squares of single letters.      \QED

        \begin{exa}
            Recall the generalized Thue-Morse word $\bt =\diamondsuit\heartsuit\triangle\diamondsuit\heartsuit\heartsuit\triangle\cdots$ that was defined on $\Sigma=\{\triangle,\diamondsuit,\heartsuit \}$ in Example \ref{E:club}. From Theorem \ref{T:critexp}, $E(\bt)=2b/m=10/3$. From Theorem \ref{T:genpos}, we compute
            $S^{3}_{5-1,3}= \{ 3,6,9,12,15,\ldots \} $
            and obtain the set of occurrences
                $A  = \{ 120, 245, 370, 495, 745,\ldots \}$
            of critical factors of length $3$.
            Here are the first two critical factors of $\bt$: \footnotesize
            $$\bt =\position{0}{\diamondsuit}\position{1}{\heartsuit}
            \position{2}{\triangle}\position{3}{\diamondsuit}
            \heartsuit\heartsuit\triangle\cdots
            \diamondsuit\heartsuit\triangle
                  \underbrace{
                  \position{120}{\triangle}\diamondsuit\heartsuit\triangle\diamondsuit
                  \heartsuit\triangle\diamondsuit\heartsuit\triangle
                  }_{\mbox{$(\triangle\diamondsuit\heartsuit)^{10/3}$}}
                  \triangle\diamondsuit\heartsuit\cdots
            \heartsuit\triangle\diamondsuit
                  \underbrace{
                  \position{245}{\diamondsuit}\heartsuit\triangle\diamondsuit\heartsuit
                  \triangle\diamondsuit\heartsuit\triangle\diamondsuit
                  }_{\mbox{$(\diamondsuit\heartsuit\triangle)^{10/3}$}}
                  \diamondsuit\heartsuit\triangle\cdots $$

        \end{exa}

\normalsize


    From Lemma~\ref{L:pos}, we obtain the following easy fact generalizing Theorem~7(a) in~\cite{jAjS00sums}, which states that $\bt$ contains the square of a single letter if and only if $\gcd(b-1,m)=1$.  

    \begin{cor}
        In the square case,
        there exists a critical factor of $\bt$
        of length $\ell$ with $b \nmid \ell$
        if and only if \emph{$\gcd(b-1,m) \mid \ell$}.
    \end{cor}

    \Proof
        If $\ell < b$, then there exists such a critical factor
        if and only if $B_\ell \neq \emptyset$, that is, if and only if
        $S^{\ell}_{b-1,m} \neq \emptyset$,
        in which case $\gcd(b-1,m) \mid \ell$.
        Otherwise, if $\ell > b$, then we know that $b=m=2$ and
        $\ell=3$, so the result follows.
    \QED

\bigskip
\noindent{\bf Acknowledgements}. The authors wish to thank the anonymous referees for carefully reading the proofs and providing helpful comments that improved the presentation. 

\end{document}